\documentclass[10pt]{article}
\usepackage{a4}
\usepackage{QED}
\usepackage{amsmath}
\usepackage{amssymb}
\usepackage{color}
\usepackage{amscd}
\usepackage{amsfonts}

\newenvironment{pro}{{\textbf{ Proof.}}}{\hfill $\blacksquare$}
\def\rd{\triangleright}
\def\red{\triangleright^*}

\def\b{\beta}

\def\G{\Gamma}
\def\D{\Delta}
\def\t{\tau}
\def\d{\delta}
\def\te{\theta}
\def\l{\lambda}

\def\m{\mu}
\def\o{\omega}

\def\s{\sigma}

\def\sou{\underline}
\def\O{\Omega}

\def\f{\rightarrow}

\def\v{\vdash}

\def\<{\langle}
\def\>{\rangle}
\def\F{\displaystyle\frac}
\def\eq{\simeq}

\newtheorem{theo}{Theorem}[section]
\newtheorem{lm}{Lemma}[section]
\newtheorem{re}{Remark}[section]
\newtheorem{cor}{Corollary}[section]
\newtheorem{df}{Definition}[section]
\newtheorem{nt}{Notation}[section]

\newtheorem{ex}{Example}[section]

\begin{document}

\vspace*{0cm}
\begin{center}
{\Large{\bf A completeness result for  the simply typed $\l \m$-calculus}}\\ [1cm]
{\bf  Karim NOUR $\&$ Khelifa SABER} \\ 
LAMA - \'Equipe LIMD \\
Universit\'e de Chamb\'ery\\
73376 Le Bourget du Lac\\
e-mail : $\{$knour@univ-savoie.fr\;\;\; ksaber@messel.emse.fr$\}$ \\[0.5cm]
\end{center}

\begin{abstract}
In this paper, we define a realizability semantics for the simply typed $\l\m$-calculus. We show that if a term is typable, then it inhabits the interpretation of its type. This result serves  to give characterizations of  the computational behavior of some closed typed terms. We also prove  a completeness result of our realizability semantics using a particular term model.
\end{abstract}

\section{Introduction}

${}\;\;\;$What came to be called the Curry-Howard correspondence has proven to be a robust technique to study proofs of intuitionistic logic, since it exhibits the structural bond between this logic and the
$\l$-calculus. T. Griffin's works \cite{4Grif} in 1990 allowed to
extend this correspondence to classical logic, which had several
consequences. On basis of this new contribution, the $\l\m$-calculus
was introduced by M. Parigot \cite{4Par1} and \cite{4Par2}. The
$\l\m$-calculus is a natural extension of the $\l$-calculus which
exactly captures the algorithmic content of proofs written in the
second order classical natural deduction system. The typed $\l
\m$-calculus enjoys all good properties: the subject reduction,
the strong normalization and confluence theorems.

The strong normalization theorem of second order classical natural
deduction \cite{4Par2} is based on a lemma known as the
correctness result, which stipulates that each term is in the
interpretation of its type.  This is also based on the notion of the
semantics of realizability. The idea of this semantics consists in
associating to each type a set of terms that {\em realizes} it, this
method has been very effective for establishing the strong
normalization of type system ``\`a la Tait and Girard''. J.- Y. Girard
used it to give a proof of the strong normalization of his system
$\cal F$, method known also as the reducibility candidates, later
M. Parigot extended this method to the classical case and provided a
proof of strong normalization of the typed $\l\m$-calculus. In a
previous work \cite{4NS1}, we adapted Parigot's method and established
a short semantical proof of the strong normalization of classical
natural deduction with disjunction as primitive.
 
In general all the known semantical proofs of strong normalization use
a variant of the reducibility candidates based on a correctness
result, which has been important also for characterizing computational behavior of some typed terms, as it was done in J.-L. Krivine's works \cite{4kri1}.
 This inspired us also to define a general semantics for classical
natural deduction in \cite{4NS2} and gave such characterizations.

The question that we now can ask is: ``does the correctness result
have a converse?''. By this we mean: ``can we find a class of types
for which the converse of the correctness result (completeness result)
holds?''.  J.R. Hindley was the first who \newpage study the completeness of
simple type systems \cite{4hind1}, \cite{4hind2} and 
\cite{4hind3}. R. Labib-sami has established in \cite{4Labib} completeness
for a class of types in Girard's system $\cal F$ known as strictely positive
types, and this for a semantics based on  sets stable under 
$\b\eta$-equivalence. S. Farkh and K. Nour revisited this
result, and generalized it, in fact they proved a refined result by
indicating that weak-head-expansion is sufficient \cite{4fn}. In
\cite{4fn2}, they established an other completeness result for a class
of types in Krivine's system $\cal AF$$2$. Recently, F. Kamareddine
and K. Nour  improved the result of Hindley, to a system with an
intersection type. Independently, T. Coquand  established in \cite{4coq}
by methods using Kripke's models, the completeness for the simply typed
$\l$-calculus.

In the present work we deal with this problem and  prove the
completeness for the simply typed $\l\m$-calculus. The semantics that
we define here is not completely different from that of \cite{4NS2} and
\cite{4NS1}, nevertheless we add a slight but an indispensable
modification to the notion of the $\m$-saturation. This semantics is inspired
 by the strong normalization proof of Parigot's $\l$$\m$-calculus, which consists in 
rewriting each reducibility candidate as a double orthogonal.  

The correcteness result allows to describe the
computational behavior of closed typed terms. We have two kinds of proofs for
such characterizations. Semantical proofs, in which we guess the computational 
behaviors, models used in such proofs are exactly built to meet the required characterization.
Syntactical proofs, where we construct the behavior based on 
the type, these proofs are shorter than the semantical ones.
 In what follows, we give at each time, both of semantics and syntactical proofs.

This paper is organized as follows. Section 2 is an introduction to the
simply typed $\l\m$-calculus. In section 3, we define the semantics and
prove its correctness. Section 4  is devoted to the
completeness result. Finally, in Section 5 we give characterizations of
some closed typed terms. 

\section {The simply typed $\l \m$-calculus}

In this work, we use the $\l \m$-calculus \`a la De Groote, where the
binder $\m$ and the naming construct are split. This allows more
expressivity than the Parigot's original version.

\begin{df} 
\begin{enumerate} 
\item Let $\mathcal{X}$ and $\mathcal{A}$ be two infinite sets of
disjoint alphabets for distinguiching $\lambda$-variables and
$\mu$-variables. The $\l \m $-terms are given by the following
grammar:

\begin{center}
$ \mathcal{T}$$ \; := $$\mathcal{X}$$ \;| \; \lambda$$\mathcal{X}$$.\mathcal{T}\;
|$$\; (\mathcal{T}\;\mathcal{T})\; |$$ \;\;\mu \mathcal{A}.\mathcal{T}$$ \;|\; (\mathcal{A}\;\mathcal{T})$
\end{center}

\item Types are formulas of the propositional logic built from the
infinite set of propositional variables $\mathcal{P} = \{X, Y,Z,...\}$
and a constant of type $\perp$, using the connective $\to$.

\item  As usual we denote by $\neg A$ the formula $A \to \perp$. 
Let $A_1,A_2,...,A_n,A$ be types, we denote the type $A_1 \to (A_2 \to
(...\to (A_n \to A)...))$ by $A_1,A_2,..., A_n\to A$.

\item Proofs are presented in  natural deduction system with two
conclusions, such that formulas in the left-hand-side of $\vdash$ are indexed by
$\l$-variables and those in right-hand-side of $\vdash$ are indexed by
$\m$-variables, except one which is indexed by a term.

\newpage

\item Let $t$ be a $\l \m $-term, $A$ a type,
$\Gamma=\{x_i:A_i\}_{1 \le i \le n}$ and $\Delta =\{a_j:B_j\}_{1 \le j
\le m}$, using the following rules, we will define ``$t$ typed with
type $ A$ in the contexts $\Gamma$ and $ \Delta$'' and we denote it $
\Gamma \vdash t:A\; ; \Delta.$

\begin{center}
 $\F{}{\Gamma\, \vdash x_i:A_i\,\, ; \, \Delta}{ax}\;\;\;$ for $1\le i
 \le n$.
\end{center}
\begin{center}
$\F{\Gamma, x:A \vdash t:B;\Delta}{\Gamma \vdash \lambda x.t:A \to B;\Delta}{\to_i}
\quad\quad\quad
\F{\Gamma \vdash u:A \to B;\Delta \quad \Gamma \vdash
v:A;\Delta}{\Gamma\vdash (u\; v):B;\Delta}{\to_e}$
\end{center}
\begin{center}
$\F{\Gamma\vdash t:\perp; \Delta, a:A}{\Gamma \vdash\mu a. t:A; \Delta}{\m}\;\;\;
\quad \F{\Gamma\vdash t:A;\Delta, a:A}{\G \v (a\; t):\perp; \D, a:A}{\perp}$
\end{center}
We denote this typed system by $S_\m$.

\item The basic reduction rules are $\beta$ and $\m$ reductions.

\begin{itemize}
\item $(\lambda x.u \; v) \triangleright_\beta u[x:=v]$
\item $(\m a.u\; v) \triangleright_\m \m a.u[a:=^*v]$

where  $u[a:=^*v]$ is obtained from $u$ by replacing inductively each
subterm in the form $(a\;w)$ in  $u$ by  $(a\;(w\;v))$.
\end{itemize}

\item We denote $t \rd t'$ if $t$ is reduced to $t'$ by one of the rules
given above. As usual $\red$ denotes the reflexive transitive closure
of $\rd$, and $\eq$ the equivalence relation induced by $\red$.
\end{enumerate}
\end{df}

We have the following results (for more details, see \cite{4Par2}).

\begin{theo}[Confluence result]
If  $t \red t_1$ and  $t \red t_2$, then there exists  $t_3$ such that  $t_1 \red t_3$ and  $t_2 \red t_3$
\end{theo}

\begin{theo}[Subject reduction]
If  $\G \vdash t : A ; \Delta$ and  $t \red t'$ then $\G \vdash t' : A ; \Delta$.
\end{theo}

\begin{theo}[Strong normalization]
If  $\G \vdash t : A ; \Delta$, then $t$ is strongly  normalizable.
\end{theo}

%\begin{re} We find in the current  literature of the $\l \m$-calculus other simplification rules such as: $\eta$, $\m'$, $\theta$, $\rho$, $\nu$,... These rules allow to get less normals forms (see \cite{4Par1}, \cite{4Par2} and \cite{4Py}).
%In this paper we do not consider these rules in order to do not complicate the proofs.

%\end{re}

\begin{df}
\begin{enumerate}
\item Let $t$ be a term and $\bar{v}$ a finite sequence of terms $($the empty sequence is denoted by $\emptyset$$)$,
then, the term $t \bar{v}$ is defined by $(t \; \emptyset) = t$ and
$(t \; u \bar{u}) = ((t \; u) \; \bar{u})$.
\item Let $t,u_1,...,u_n$ be terms and $\bar{v}_1,...,\bar{v}_m$
finite sequences of terms, then 

$t[(x_i:=u_i)_{1 \le i\le n};(a_j:=^*\bar{v}_j)_{1 \le j \le m}]$ is
obtained from the term $t$ by replacing inductively each $x_i$ by
$u_i$ and each subterm in the form $(a_j\;u)$ in $t$ by
$(a_j\,(u\,\bar{v}_j))$.
\end{enumerate}
\end{df}

\begin{re} In order to avoid the heavy notation of the substitution\\ $[(x_i:=u_i)_{1 \le i\le n};(a_j:=^*\bar{v}_j)_{1 \le j \le m}]$, we denote it by $\sigma$ $($which is not an object of the syntax$)$. Then $t[(x_i:=u_i)_{1 \le i\le n};(a_j:=^*\bar{v}_j)_{1 \le j \le m}]$ is denoted by $t\sigma$.

\end{re}

\begin{lm}\label{sigma} 
Let $t$, $t'$ be  terms and $\sigma$ a substitution, if $t\red t'$, then, $t\sigma \red t'\sigma$.
\end{lm}

\begin{pro} By induction on $t$.\end{pro}

%%%%%%%%  la semantik%%%%%%%%%%%%%%%%%%

\section{The  semantics of $S_\m$}

In this part we define the realizability semantics and prove its
correctness.

\begin{df}
\begin{enumerate}
\item We say that a set of terms ${\cal{S}}$ is saturated when the conditions: $v \red u$ and $ u\in {\cal{S}}$ imply   $v\in {\cal{S}}$ for all terms $u$ and $v$.

\item Let us take a saturated set of terms ${\cal{S}}$ and a set
$\mathcal{C}$ of an infinite classical variables
$($$\m$-variables$)$. We say that ${\cal{S}}$ is
$\mathcal{C}$-saturated when the condition: $t \in {\cal{S}}$ implies
$\mu a.t \in {\cal{S}}$ and $(a\; t) \in {\cal{S}}$ for all term $t$
and all $\m$-variable $a\in \mathcal{C}$
\end{enumerate}
\end{df}

\begin{re}\label{satur} 
The difference between this semantics and those defined in \cite{4NS2}
and \cite{4NS1}, is the notion of the $\cal{C}$-saturation which is not
necessary for the correctness part, but indispensable for the
completeness side.  It is obvious that this notion introduces
ill-typed terms, thing which seems to go against
completeness. Nevertheless, the key point is that $\cal{C}$ is a
parameter attached to a particular model, therefore when we take the
intersection of all models, all these bad terms are removed. This is
exaclty what is done in the proof of the theorem \ref{comp}.
\end{re}

\begin{df}
\begin{enumerate}

\item Consider two sets of terms $\cal{K}$ and $\cal{L}$, we define a
new set of terms: ${\cal{K}} \leadsto {\cal{L}} =\{ t$ / $(t\; u) \in
{\cal{L}},$ for each $ u \in \cal{K}\}$. It is clear that when
$\cal{L}$ is a saturated set, then $\cal{K}\leadsto \cal{L}$ is also 
saturated one.

\item We denote $\mathcal{T} \cup \mathcal{A}$ by $\mathcal{T'}$ and
$\mathcal{T'}^{<\omega}$ the set of finite sequences of elements of
$\mathcal{T'}$. Let $t$ be a term and $\pi \in \mathcal{T'}^{<\o}$,
then the term $(t\;\pi)$ is defined by $(t\;\emptyset) = t$, $(t\;\pi)
= ((t\;u)\; \pi')$ if $\pi= u\pi'$ and $(t\;\pi) = ((a\;t)\;\pi')$ if
$\pi= a\pi'$.

\item Let ${\cal{S}}$ be a set of terms and $\cal{X}$$ \subseteq
 \cal{T'}^{<\omega}$, then we define $\cal{X}$$ \leadsto {\cal{S}}
 =\{t\; /\; (t\; \pi) \in {\cal{S}},$ for each $\pi \in \cal{X}\}$.
\end{enumerate}
\end{df}

\begin{re} 
The fact that the application $(a\,t)$ is denoted by $(t\,a)$ is not
something new, it is already present in Saurin's work
\cite{4Sau}. Except that for us, it is a simple notation in order to
uniformize the definition of the application. But for Saurin, it is
crucial to obtain the separation theorem in the $\l\m$-calculus.
\end{re}

\begin{df}\label{modl} 
Let ${\cal{S}}$ be a $\mathcal{C}$-saturated set and
$\{{\cal{R}}_i\}_{i \in I}$ subsets of terms such that ${\cal{R}}_i =
{\cal{X}}_{{\cal{R}}_i} \leadsto {\cal{S}}$ for some
${\cal{X}}_{{\cal{R}}_i} \subseteq \mathcal{T'}^{<\omega}$. A model
$\mathcal{M}$=$\<\mathcal{C}$$,{\cal{S}},\{{\cal{R}}_i\}_{ i \in I}\>$
is the smallest set containing ${\cal{S}}$ and ${\cal{R}}_i$, and
closed under the constructor $\leadsto$.
\end{df}

\begin{lm} 
Let $\mathcal{M} = \<\mathcal{C}$$,{\cal{S}},\{{\cal{R}}_i\}_{i \in
I}\>$ be a model and ${\cal{G}} \in \mathcal{M}$. There exists a set
$\cal{X_G} \subseteq \mathcal{T'}^{<\o} $ such that ${\cal{G} = X_G}
\leadsto {\cal{S}}$.
\end{lm}

\begin{pro} 
By induction on $\cal{G}$.
\begin{itemize}
\item[-] If $\cal{G}={\cal{S}}$, take  $\cal{X_G} = \{\phi\}$.
\item[-] If ${\cal{G}}={\cal{R}}_i$, take ${\cal{X_G}} = {\cal{X}}_{{\cal{R}}_i}$.
\item[-] If ${\cal{G}}={\cal{G}}_1 \leadsto {\cal{G}}_2$, then, by induction
hypothesis, ${\cal{G}}_2={\cal{X}}_{{\cal{G}}_2}\leadsto {\cal{S}}$
where ${\cal{X}}_{{\cal{G}}_2} \subseteq \mathcal{T'}^{<\o}$, and
take ${\cal{X_G}}=\{u \bar{v}$ / $u\in {\cal{G}}_1$ and $\bar{v}\in
{\cal{X}}_{{\cal{G}}_2}\}$.
\end{itemize}
\end{pro}

\newpage
\begin{df} 
Let $\mathcal{M} = \<\mathcal{C},{\cal{S}},\{{\cal{R}}_i\}_{i \in
I}\>$ be a model and ${\cal{G}} \in \mathcal{M}$. We
define the set ${\cal{G}}^\perp =\cup\{{\cal{X_G}}$ / ${\cal{G=
X_G}}\leadsto {\cal{S}} \}$.
\end{df}

\begin{lm}\label{orth} 
Let $\mathcal{M} = \<\mathcal{C},{\cal{S}},\{{\cal{R}}_i\}_{i \in
I}\>$ be a model and ${\cal{G}} \in \mathcal{M}$. We have
${\cal{G = G}} ^\perp \leadsto {\cal{S}}$.
\end{lm}

\begin{pro} Immediate.
%This comes from the fact that: if for every $j\in J$, ${\cal{G}}={\cal{X}}_{{\cal{G}}_j} \leadsto {\cal{S}}$,
%then, ${\cal{G}} = ({\cup}_{j\in J} {\cal{X}}_{{\cal{G}}_j}) \leadsto
%{\cal{S}}$.
\end{pro}

\begin{df} 
\begin{enumerate}
\item Let $\mathcal{M} =\<\mathcal{C},{\cal{S}},\{{\cal{R}}_i\}_{i \in
I}\>$ be a model. An $\mathcal{M}$-interpretation $\cal{I}$ is an
application $X \mapsto {\cal{I}}(X)$ from the set of propositional
variables $\mathcal{P}$ in $\mathcal{M}$ which we extend for any
formula as follows:
\begin{itemize}
\item ${\cal{I}}(\perp)={\cal{S}}$
\item ${\cal{I}}(A \to B)= {\cal{I}}(A) \leadsto {\cal{I}}(B)$.
\end{itemize}
\item For any type $A$, we denote $\vert A \vert_\mathcal{M}$$ =\bigcap
\{ {\cal{I}}(A)$ / ${\cal{I}}$ an $\mathcal{M}$-interpretation$\}$.
\item For any type $A$, $|A|= \bigcap\{\vert A \vert_\mathcal{M}$ / $\mathcal{M}$ a model$\}$.
\end{enumerate}
\end{df}

The notion of $\mathcal{C}$-saturation is indispensable for
completeness but, as we said in the remark \ref{satur}, it provides
ill-terms. The presence of such terms has some drawbacks on the
correctness side, hence we introduce in the following definition a
parameterized relation $\hookrightarrow_{\mathcal{C}}$.

\begin{df}

Let $u,\, v$ be two terms. The expression $u
\hookrightarrow_{\mathcal{C}} v$ means that $v$ is obtained from $u$
by replacing the free classical variables of $u$ by some others in
$\mathcal{C}$, i.e, if we denote $u$ by $u[a_1,...,a_n]$ where the
$a_i$ are the free classical variables of $u$, then $v$ will be
$u[a_1:=b_1,...,a_n:=b_n]$ where $b_i \ne b_j$ for $(i \ne j)$ and
$b_i \in \mathcal{C}$ for each $1 \le i \le n$ $($it is obvious that
$\hookrightarrow_{\mathcal{C}}$ is parameterized by $\mathcal{C}$$)$.

\end{df}

\begin{lm} [Correctness] \label{adq} 
Let $\Gamma =\{x_i : A_i\}_{\substack {1\le i\le n}}$, $\Delta =\{a_j
: B_j\}_{\substack {1\le j\le m}}$,
$\mathcal{M}=\<\mathcal{C},{\cal{S}},\{{\cal{R}}_i\}_{i \in I}\>$ a
model, ${\cal{I}}$ an $\mathcal{M}$-interpretation, $u_i \in
\;{\cal{I}}(A_i)$, $\bar{v_j} \in\; ({\cal{I}}(B_j))^\perp$, $\sigma
=[(x_i:=u_i)_{1 \le i \le n};(a_j:=^*\bar{v}_j)_{1\le j\le m}]$, and
$u,v$ two terms such that $u \hookrightarrow_{\mathcal{C}} v$. If
$\Gamma \vdash u:A \,\,\,;\,\Delta$, then, $v\sigma \in {\cal{I}}(A)$.
\end{lm}

\begin{pro} 
By induction on the  derivation, we consider the last used rule.

\begin{enumerate}

\item[$ax$:] In this case $u=x_i=v$ and $A=A_i$, then $v\s =u_i \in {\cal{I}}(A)$.

\item[$\to_i$:] In this case $u= \l x.u_1$ and $A =B \to C$ such that
$\G ,x:B \v u_1:C \; ;\D$. Then $v =\l x.v_1$ and $u_1
\hookrightarrow_{\mathcal{C}} v_1$. Let $w \in {\cal{I}}(B)$ and $\d
=\s+[x:=w]$, by induction hypothesis, $v_1\d \in {\cal{I}}(C)$, hence
$(\l x.v_1\s \; w) \in {\cal{I}}(C)$, therefore $\l x.v_1 \s \in
{\cal{I}}(B) \leadsto {\cal{I}}(C)$. Finally $v\s \in {\cal{I}}(A)$.

\item[$\to_e$:] In this case $u =(u_1\;u_2)$, $\G \v u_1:B \to A \;
  ;\D$ and $\G \v u_2:B\; ;\D$. We also have  $v=(v_1\;v_2)$ where $u_1
  \hookrightarrow_{\mathcal{C}}v_1$ and $u_2
  \hookrightarrow_{\mathcal{C}} v_2$. By induction hypothesis, $v_1\s
  \in {\cal{I}}(B) \leadsto {\cal{I}}(A)$ and $v_2\s \in
      {\cal{I}}(B)$, therefore $(v_1\s\;v_2\s) \in {\cal{I}}(A)$, this
      implies that $v\s \in {\cal{I}}(A)$.

\item[$\m$:] In this case $u =\m a. u_1$, then $v=\m b.v_1$
where $u_1 \hookrightarrow_{\mathcal{C}} v_1$ and $b$ is a new variable which
belongs to $\mathcal{C}$ and not free in $u_1$ (there is always such
variable because $\mathcal{C}$ is infinite). Let $\bar{v}\in
({\cal{I}}(A))^{\perp}$ and $\d =\s +[b:=^*\bar{v}]$. By  induction hypothesis,
$v_1\d \in {\cal{S}}$, and by the definition of
${\cal{S}}$, we have, $\m b. v_1\d \in {\cal{S}} $. Since
$(\m b. v_1 \s \; \bar{v}) \red \m b. v_1 \d$, then, $\m
b. v_1 \s \in {\cal{I}}(A)$, i.e, $v\s \in {\cal{I}}(A)$.

\item[$\bot$:] In this case $u =(a\; u_1)$, then, $v=(b\;v_1)$ where
$u_1 \hookrightarrow_{\mathcal{C}} v_1$ such that the free variable
$a$ was replaced by $b$ in $u_1$ and $b \notin Fv(u_1)$ is new
variable which belongs to $\mathcal{C}$. Let $\d =\s+[b:=^*\bar{v}]$
where $\bar{v} \in ({\cal{I}}(A))^{\perp}$, by induction hypothesis,
$v_1\d \in {\cal{I}}(A)$, hence $(v_1\d \;\bar{v}) \in
{\cal{S}}$. Therefore, by the definition of ${\cal{S}}$, $(b\;(v_1\d
\;\bar{v})) \in {\cal{S}}$, finally $v\s \in {\cal{S}}$.

\end{enumerate}
\end{pro}

\begin{cor} \label{closed}
Let $A$ be a type and $t$ a closed term. If $\vdash t:A$, then, $t \in \vert A \vert$.
\end{cor}

\begin{pro} 
Let $\mathcal{M}$ be a model and ${\cal{I}}$ an
$\mathcal{M}$-interpretation. Since $\vdash t:A$, then, by the lemma
\ref{adq}, $t \in {\cal{I}}(A)$. This is true for any model
$\mathcal{M}$ and for any $\mathcal{M}$-interpretation $I$, therefore
$t \in |A|$.
\end{pro}

%%%%%%%%%%%%%%%%%typage obstacle%%%%%%%%%%%%%%%%%%%%%%

\section{The completeness result}

Roughly speaking, completeness of the semantics amounts to saying that
if $t$ is in the interpretation of a type $A$, then $t$ has the type
$A$. In order to prove the completeness result, we construct in the
following part a particular term model.

\begin{df}$($and notation$)$\label{compl}
\begin{enumerate}
\item Let $\Omega =\{x_i$ / $i\in \mathbb{N}\}$ $\cup$ $\{a_j$ /
$j\in \mathbb{N}\}$ be an enumeration of infinite sets of $\l$ and
$\mu$-variables.
\item Let $\O_1=\{A_i$ / $i\in \mathbb{N}\}$ be an enumeration of all
  types where each type comes infinitely many times.

\item Let $\O_2= \{B_j$ / $j\in \mathbb{N}\}$ be an enumeration of all
types where the type $\perp$ comes  infinitely  many times.

\item We define $\mathbb{G} =\{x_i : A_i$ / $i\in \mathbb{N}\} $ and $\mathbb{D}
=\{a_j : B_j$ / $j\in \mathbb{N}\}$.

\item Let $u$ be a  term, such that $Fv(u) \subseteq $ $ \Omega$, the
 contexts $\mathbb{G}_u$ $($resp $ \mathbb{D}_u$$)$ are defined as the
 restrictions of $\mathbb{G}$ $($resp $\mathbb{D}$$)$ at the declarations
 containing the variables of  $Fv(u)$.

\item The notation $\mathbb{G} \vdash u:C ;\,\,\mathbb{D}$ means that
$\mathbb{G}_u \vdash u:C ;\,\,\mathbb{D}_u $, we denote $\mathbb{G}
\vdash^* u:C ;\,\,\mathbb{D}$ iff there exists a term $u'$, such that
$u \red u'$ and $\mathbb{G} \vdash u':C ;\,\,\mathbb{D}$.

\item Let $\mathbb{C} =\{a_j$ / $(a_j : \,\perp) \in \mathbb{D}\}$ and
$\mathbb{S}=\{ t$ / $\mathbb{G} \vdash^* t: \perp ;\,\,\mathbb{D}\}$. 

\item For each propositional variable $X$, we define a set of terms
$\mathbb{R}_X =\{t$ / $\mathbb{G} \vdash^* t: X ;\,\,\mathbb{D}\}$.
\end{enumerate}
\end{df}

\begin{lm}
\begin{enumerate}
\item $\mathbb{S}$ is a $\mathbb{C}$-saturated set.
\item The sets $\mathbb{R}_X$ are saturated.
\item For each propositional variable $X$, $\mathbb{R}_X =\{ a_j$ /
$(a_j: X) \in {\mathbb{D}} \} \leadsto \mathbb{S}$.
\item  $\mathbb{M} =\<\mathbb{C},\mathbb{S},{(\mathbb{R}_X)}_{X \in \mathcal{P}}\>$ is a model
\end{enumerate}
\end{lm}

\begin{pro} Easy.
\end{pro}
\begin{re} 
Observe that the model $\mathbb{M}$ is parameterized by the two
infinite sets of variables and the enumerations, we need just these
infinite sets of variables and not all the variables. This is an
important remark since it will serve us in the proof of the theorem
\ref{comp}.
\end{re}

\begin{df}
We define the $\mathbb{M}$-interpretation $\mathbb{I}$ as follows:
\begin{itemize}
\item $\mathbb{I}(\perp)= \mathbb{S}$.
\item $\mathbb{I}(X)= \mathbb{R}_X$ for each propositional variable.
\end{itemize}
\end{df}

\begin{lm} 
Let $y$ be a $\l$-variable, $\s = [(x_i:=y)_{1 \leq i \leq n},
(a_i:=^* y)_{1 \leq j \leq m}]$ a substitution and $t$ a term.
\begin{enumerate}
\item If $(t\s\; y)$ is normalizable, then $t$ is normalizable.
\item If $t\s$ is normalizable, then $t$ is normalizable.
\end{enumerate}
\end{lm}

\begin{pro} By a simultaneous induction on $t$, we use  the
standardization theorem of the $\l\m$-calculus \cite{4Py}.

\begin{enumerate}
\item We examine the case where $t= \l x. u$. Then $(t\s\;y) = (\l
x. u\s\;y)$ is normalizable, this implies that $u\s[x:=y]$ is
normalizable, hence by $(2)$, $u$ is normalizable, therefore $t$ is
normalizable too.

\item We examine the case where $t = (a\, u)$. Then $t\s
 =(a\;(u\s\;y))$ is normalizable, this implies that $(u\s\;y)$ is
 normalizable, hence by $(1)$, $u$ is normalizable, therefore $t$ is
 normalizable too.
\end{enumerate}
\end{pro}

\begin{cor}\label{norm} 
Let $t$ by a term and $y$ a $\l$-variable.
If $(t\;y)$ is normalizable, then, $t$ is also normalizable.
\end{cor}

\begin{pro} 
Immediate from the previous lemma.
\end{pro}

\begin{lm} \label{perdre patience} 
Let $t$ and $\t$ be two normal terms, $y$ a $\l$-variable such that $y
\notin Fv(t), \;\; (t\; y) \red \t$, $A$ and $B$ types, and $\Gamma ,
y:A \v \t: B; \; \Delta$. Then $\Gamma\, \v t: A \to B; \;
\Delta$.
\end{lm}

\begin{pro} 
See the appendix.
\end{pro}

\begin{lm}\label{fatiguant}
Let $A$ be a type and $t$ a term.
\begin{enumerate}
\item If $\mathbb{G} \, \vdash^* t:A \; ;\mathbb{D}$, then $t \in \mathbb{I}(A)$.
\item If $t \in \mathbb{I}(A)$, then $\mathbb{G} \, \vdash^* t:A \; ; \mathbb{D}$.
\end{enumerate}
\end{lm}

\begin{pro} By  a simultaneous induction on the type $A$.\\

{\underline {Proof of $(1)$}}

\begin{enumerate}
\item If $A=X$ or $\perp$, the result is immediate from the definition of $\mathbb{I}$.
\item Let $A=B \to C$ and $\mathbb{G} \, \vdash^* t:A \; ; \mathbb{D}$, then $t \red
  t'$ such that: $\mathbb{G} \, \vdash t': B \to C \; ;
  \mathbb{D}$. Let $u \in \mathbb{I}(B)$. By  induction
  hypothesis $(2)$, we have $\mathbb{G} \, \vdash^* u:B \; ;
  \mathbb{D}$, this implies that $u \red u'$ and $\mathbb{G} \, \vdash
  u': B \; ; \mathbb{D}$. Hence $\mathbb{G} \, \vdash (t'\;u'): C \; ;
  \mathbb{D}$, so, by the fact that $(t\;u) \red (t'\;u')$, we
  have $\mathbb{G} \, \vdash^* (t\;u):C \; ; \mathbb{D}$, then, by
  induction hypothesis $(1)$, $(t\;u) \in \mathbb{I}(C)$. Therefore
  $t \in \mathbb{I}(B \to C)$.
\end{enumerate}

\underline{Proof of $(2)$}

\begin{enumerate}
\item If $A=X$ or $\perp$, the result is immediate from the definition of $\mathbb{I}$.
\item Let $A=B \to C$, $t \in \mathbb{I}(B) \leadsto \mathbb{I}(C)$ and $y$ be a $\l$-
  variable such $y \not\in Fv(t)$ and $(y:B) \in \mathbb{G}$. We have $y:B
\v y:B$, hence, by induction hypothesis
  $(1)$, $y \in \mathbb{I}(B)$, then, $(t\; y) \in
  \mathbb{I}(C)$. By induction hypothesis $(2)$, $\mathbb{G} \,
  \vdash^* (t\;y) :C \; ; \mathbb{D}$, then $(t\; y) \red t'$ such
  that $\mathbb{G} \, \vdash t':C \; ; \mathbb{D}$ and, by the
  corollary \ref{norm}, $t$ is a normalizable term. The normal form of
  $t$ can be either $(x\; u_1)\;u_2...u_n$ either $\l x.u$ or $\m a.u$
  $($the case $(a \;u)$ gives a contradiction for  typing reasons$)$.

\begin{enumerate}
\item If $t \red (x\; u_1)\;u_2...u_n$ with $u_i$ normal terms, then
 $\mathbb{G} \, \vdash (x\;u_1)\;u_2...u_n y:C \; ; \mathbb{D}$, $x
 : E_1,E_2,...,E_n \to (B \to C) \in \mathbb{G}$, $\mathbb{G} \,
 \vdash u_i:E_i \; ; \mathbb{D}$ and $\mathbb{G} \, \vdash y:B \; ;
 \mathbb{D}$.  Therefore $\mathbb{G} \, \vdash (x\; u_1)\;u_2...u_n :B
 \to C \; ; \mathbb{D}$, finally $\mathbb{G} \, \vdash^* t :B \to
 C \; ; \mathbb{D}$.

\item If $ t \,\red\, \l x.u$ where $u$ is a normal
  term, then, since $\mathbb{G}$ contains an infinite number of declarations
  for each type, let $y$ be a $\l$-variable such that
  $(y:B) \in \mathbb{G}$ and $y \notin Fv(u)$. We have $(t\;y)
  \red u[x:=y]$ and $\mathbb{G} \, \vdash u[x:=y]:C  \; ;
  \mathbb{D}$, hence $\mathbb{G} \, \vdash \l y.u[x:=y]:B \to C  \; ;
  \mathbb{D}$ and, by the fact that $y \notin Fv(u)$, 
  $\l y.u[x:=y] = \l x.u $. Therefore $\mathbb{G} \, \vdash \l x.u:B \to C
\; ; \mathbb{D}$, finally $\mathbb{G} \, \vdash^* t:B \to C  \; ; \mathbb{D}$.

 \item If $t \,\red \,\m a.u$ where $u$ is a normal term, then let $y$ be
 a $\l$-variable such that $(y:B) \in \mathbb{G}$ and $y \notin
 Fv(u)$. We have $(t\;y) \red \mu a.u[a:=^*y] \red \mu a.u'$ where
 $u'$ is the normal form of $u[a:=^*y] $, so we have $\mathbb{G} \,, y:B
 \vdash \mu a.u':C \; ; \mathbb{D}$. By the lemma \ref{perdre
 patience}, we obtain $\mathbb{G} \, \vdash \mu a.u: B \to C \; ;
 \mathbb{D}$, finally $\mathbb{G} \, \vdash^* t: B \to C \; ;
 \mathbb{D}$.
\end{enumerate}

\end{enumerate}
\end{pro}

\begin{theo} \label{comp} 
Let $A$ be a type and $t$ a term. We have $t \in \vert A \vert$ iff
there exists a closed term $t'$ such that $t\red t'$ and $\v t':A$.
\end{theo}

\begin{pro}
$\Leftarrow)$ By the lemma \ref{adq}.

$\Rightarrow)$ We consider an infinite set of $\l$ and $\m$ variables $\Omega$ such
that it contains none of the free variables of $t$, then from this set
we build the completeness model as described in the definition
\ref{compl}. If $t\in \vert A \vert$, then $t \in \mathbb{I}(A)$,
hence by $(1)$ of the lemma \ref{fatiguant} and by the fact that
$Fv(t') \subseteq Fv(t)$, we have $t \red t'$ and $\vdash t':A$.
\end{pro}

\begin{cor}
 Let $A$  be a  type and $t$ a term.
\begin{enumerate}
\item If $t \in \vert A \vert$, then $t$ is normalizable.
\item If $t \in \vert A \vert$, then there exists a closed term $t'$
such that $t \eq t'$.
\item $\vert A \vert$ is closed under equivalence. 
\end{enumerate}
\end{cor}

\begin{pro} 
$(1)$ and $(2)$ are  direct consequences of theorem \ref{comp}.  $(3)$
can be deduced from the theorem \ref{comp} and the lemma \ref{adq}.
\end{pro}

\section{Characterization of some typed terms}

We begin by adding to our system new propositional
constants to obtain a new parameterized typed system.  In such systems
we can characterize the syntactical form of a term having some type,
this will be useful for the proof of the lemma \ref {obstacle}. This
part is inspired by Nour's works \cite{4Nour5} and \cite{4Nour6}.

\subsection{The system ${S_\m}^{\bar{O}}$}

\begin{df} Let $\bar{O}=O_1,...,O_n$ be a sequence of
fresh propositional constants.
\begin{enumerate}
%\item We said that $\bar{O}$ is different from $\perp$ iff each $O_i$
%is different from $\perp$. 
\item A type $A$ is said an  $\bar{O}$-type iff $A$ is obtained by the following rules:
\begin{itemize}
\item Each $O_i$ is an $\bar{O}$-type.
\item If $B$ is an $\bar{O}$-type, then, $ A \to B$ is an $\bar{O}$-type.
\end{itemize}

\item The typed system ${S_\m}^{\bar{O}}$ is the system
${S_\m}$ at which we add the following conditions:
\begin{itemize}

\item The rules $ax$ is replaced by

\begin{center}
 $\F{}{\Gamma\, \vdash_{\bar{O}} x_i:A_i\,\, ; \, \Delta}{ax}$
\end{center}
where $\D$ does not contain declarations of the form $a:C$ such that $C$ is an $\bar{O}$-type.
\item The rules $\to_e$ is replaced by
\begin{center}
$\F{\Gamma \v_{\bar{O}}  u:A \to B;\Delta\;\;\;\;\Gamma \v_{\bar{O}}
v:A;\Delta}{\Gamma \v_{\bar{O}} (u\;v):B ;\Delta}{\to_e}$ 
\end{center}
where $B$ is not an $\bar{O}$-type.
\end{itemize}

\end{enumerate}
\end{df}

\begin{re}
It is obvious that ${S_\m}$$^{\bar{O}}$ can be seen as the system
${S_\m}$ where the syntax of formulas is extended by the new constants
$\bar{O}$ and some restrictions are imposed on the typing
rules. Therefore in the remainder of this work we consider that, any
typed term in the system ${S_\m}$$^{\bar{O}}$ is strongly
normalizable.
\end{re}

\begin{lm}\label{toto} 
If $\G \v t: A\;; \D$, $X$ a propositional variable and $F$ is not an
$\bar{O}$-type, then $\G \v_{\bar{O}} t: A[X:=F]\;; \D$.
\end{lm}

\begin{pro} 
By induction on the derivation.
\end{pro}\\

The following lemma stipulates that the new system
${S_\m}$$^{\bar{O}}$ is closed under reduction (subject reduction).

\begin{lm} 
If $\G \v_{\bar{O}} t: A\;; \D$ and $t \red t'$, then $\G \v_{\bar{O}}
t': A\;; \D$
\end{lm}

\begin{pro} 
By induction on the length of the reduction $t \red t'$. It suffices
to check this result for $t \rd_\b t'$ and $t \rd_\m t'$. We process
by induction on $t$.
\end{pro}

\begin{lm}\label{obstacle}
Let $\G=\{x_i:A_i\}_{1\le i \le n}$, $\D=\{a_j:B_j\}_{1\le j \le m}$
$\bar{O}=O_1,...,O_k$ and $1 \leq l \leq k$. If $\G \v_{\bar{O}} t:O_l
\; ;\D$, then, $t=x_j$ for some $1\le j\le n$ and $A_j= O_l$.
\end{lm}

\begin{pro} 
By induction on the derivation.
\begin{itemize}
\item[$ax$:] Then, $\G \v x_j:A_j;\; \D$, hence $t=x_j$ and $O_l = A_j$.
\item[$\to_i$:] A contradiction because this implies that $O_l$ is not
  atomic.
\item[$\to_e$:] This implies that $t=(u\;v)$, then, $\G \v u: A \to
O_l;\D$, therefore this gives a contradiction with the restriction on
the rule $\to_e$ since $O_l$ is an ${\bar{O}}$-type.
\item[$\m$:] Then, $t= \m a.t_1$ and $\G \v t_1:\perp;\; \D', a:O_l$,
where $\D = \D' \cup \{a:O_l\}$, therefore this gives a contradiction
with the fact that $\D$ does not contain declarations of the form
$a_j:O_j$.
\item[$\perp$:] A contradiction because $O_l$ is different from $\perp$.
\end{itemize}
\end{pro}\\

Now we give some applications of the lemma \ref{adq}. We will see that
the operational behavior of a typed term depends in ``certain sense'' only of its type.

\begin{df}
Let $t$ be a term. We denote $M_t$ the smallest set containing $t$
such that: if $u \in M_t$ and $a \in {\cal A}$, then $\m a.u \in
M_t$ and $(a \; u) \in M_t$. Each element of $M_t$ is denoted
$\sou{\m}.t$. For example, the term $\m a.\m b.(a \; (b \; (\m c.
(a \; \m d.t))))$ is denoted by $\sou{\m}.t$.
\end{df}

\subsection{Terms of type $\perp \to X$}

\begin{ex}
Let $e_1 =\lambda x.\mu a.x$ and $e_2 = \lambda x.\mu b.(b\, \mu
a.x)$, we have:\\ $\vdash e_i : \perp \to X$.

Given a $\l$-variable $x$, and  a finite sequence of $\l$-variables $\bar{y}$, we have:
\begin{itemize}
\item $(e_1\; x) \;\bar{y} \red \m a. x$
\item $(e_2\; x) \;\bar{y} \red \m b.(b \; \mu a. x)$
\end{itemize}
\end{ex}

The operational behavior of closed terms with the type $\perp \to X$
is given in the following theorem.

\begin{theo}
Let $e$ be a closed term of type $\perp \to X$, then, for each
$\l$-variable $x$ and for each finite sequence of $\l$-variables
$\bar{y}$, $(e \;x)\; \bar{y} \red \sou {\m}. x$
\end{theo}

\begin{pro}

\underline{Semantical proof:}

Let $x$ be a $\l$-variable and $\bar{y}$ a finite sequence of
$\l$-variables. Let $\mathcal{C}=\mathcal{A}$, take
${\cal{S}}=\{t\,/\, t\red \sou{\m}. x \}$ and
${\cal{R}}=\{\bar{y}\}\leadsto {\cal{S}}$. It is clear that
${\cal{S}}$ is $\mathcal{C}$-saturated set and $x \in {\cal{S}}$. So
let $\mathcal{M}=\<\mathcal{C},{\cal{S}},{\cal{R}}\>$ and take
${\cal{I}}$ the interpretation which at $X$ associates
${\cal{I}}(X)={\cal{R}}$. By the lemma \ref{adq}, $e \in {\cal{I}}
(\perp \to X)$, then, $e\in {\cal{S}}\leadsto {\cal{R}}$, i.e, $e \in
{\cal{S}} \leadsto (\{\bar{y}\}\leadsto {\cal{S}})$, therefore $(e\;x)
\in \{\bar{y}\} \leadsto {\cal{S}}$, and $(e\;x)\;\bar{y} \in
{\cal{S}}$. Finally $(e \;x)\; \bar{y} \red \sou {\m}. x$.\\

\underline{Syntactical proof:}

We can also give a syntactical proof of this result. Let
$\bar{O}=O_1,...,O_n$ be a sequence of new constants, $A=O_1,...,O_n
\to \bot$ and $\bar{y} = y_1...y_n$ a sequence of $\l$-variables.  By
the lemma \ref{toto}, $\v_{\bar{O}} e: \perp \to A$, then, $x:\perp,
(y_i:O_i)_{1\le i\le n} \v_{\bar{O}} (e\;x)\bar{y}:\perp$, hence
$(e\;x)\bar{y}\, \red \t$. It suffices to prove that, if $\t$ is a
normal term and $x:\bot, (y_i:O_i)_{1\le i\le n} \v_{\bar{O}}
\t:\bot\; ; (b_j: \bot)_{1 \le j \le m}$, then $\t =\sou{\m}.x$. This
can be proved easily by induction on $\t$.
\end{pro}

\begin{cor}
Let $e$ be a closed term of type $(\perp \to X)$, then, for each term
$u$ and for each $\bar{v} \in \mathcal{T}^{<\o}$, $(e \,u)\, \bar{v}
\red \sou {\m}.u$
\end{cor}

\begin{pro} 
Immediate from the previous theorem and the lemma \ref{sigma}.
\end{pro}

\begin{re} 
Let $\v e: \bot \to X$, the term $(e \;u)$ modelizes an instruction
like ${\tt exit}(u)$ $(${\tt exit} is to be understood as in the {\tt C}
programming language$)$. In the reduction of a term, if the subterm
$(e \;u)$ appears in head position $($the term has the form $($$(e\; u) \;
\bar{v})$$)$, then, after some reductions, the sequence $\bar{v}$ is
deleted, and we obtain $\sou{\m}.u$ as result.
\end{re}

\subsection{Terms of  type $(\neg X \to X) \to X$}

\begin{ex}  
Let the terms $E_1 =\l x.\m a.(a\,(x\,\, \l z.(a\; z)))$\\ and $E_2
=\l x.\mu a.(a\,\,(x\,\,(\l z_1.(a (x\,\, \l z_2. (a\,z_1)))))) $,\\
we have: $\,\,\vdash E_i : (\neg X \to X) \to X$.\\ Given $\l$-variables
$x, z_1, z_2$ and a finite sequence of $\l$-variables $\bar{y}$, we
have:

\begin{itemize}
\item $(E_1\; x)\; \bar{y} \red \m a.(a\;((x\; \te_1)\;\bar{y}))$ and
$(\te_1\; z_1) \red (a\,(z_1\;\bar{y}))$, where $\te_1 = \l
z.(a\,(z\,\bar{y}))$.
\item $(E_2\; x) \;\bar{y} \red \mu a.(a \,((x\;\te_1)\;\bar{y}))$,
$(\te_1\,\, z_1) \red (a\; ((x \;\te_2)\;\bar{y}))$, and $(\te_2
\;z_2) \red (a\;(z_1 \;\bar{y}))$, where $\te_1 = \l z_1.(a\,((x\, \l
z_2.(a\,(z_1\,\bar{y})))\,\bar{y}))$ and $\te_2 = \l z_2.(a\,(z_1
\,\bar{y}))$.
\end{itemize}
\end{ex}

The following theorem describes the computational behavior of closed
terms with type $(\neg X \to X) \to X$.

\begin{theo}
Let $E$ be a closed term of type $(\neg X \to X) \to X$, then, for
each $\l$-variable $x$, for each finite sequence of $\l$-variables
$\bar{y}$ and for each sequence of $\l$-variables $(z_i)_{i \in
  {\mathbb{N}^*}}$ such that: $x$, $y_j$ are differents from any
$z_i$. There exist $m \in {\mathbb{N}}^*$ and terms $\te_1,...,\te_m$,
such that we have:
\begin{itemize}
\item$(E\; x) \bar{y} \red \sou{\mu}. (x\,\, \te_1)\;\bar{y} $
\item$(\te_k\; z_k)\red \sou{\mu}. (x\,\, \te_{k+1})\;\bar{y}$ $\;\;\;\;$
  for all  $1\le k \le m-1 $
\item$(\te_m\; z_m) \red \sou{\mu}. (z_l\,\, \bar{y})$
  $\;\;\;\;\;\;\;$ for some  $1\le l \le m $
\end{itemize}
\end{theo}

\begin{pro}

\underline{Semantical proof:}

Let $x$ be a $\l$-variable, $\bar{y}$ a finite sequence of
$\l$-variables and $(z_i)_{i \in \mathbb{N}^*}$ a sequence of
$\l$-variables as in the theorem above. Take ${\cal{S}}=\{t$
/$\forall\, r\ge 0$: Either $[\exists m\ge 1,\; \exists
\te_1,...,\te_m, \;\exists\, j$: $t \,\red \sou{\m}. ((x
\;\te_1)\,\bar{y})$, $(\te_k\; z_{k+r})\, \red\, \sou{\m}.  ((x
\,\te_{k+1})\,\bar{y})$ for every $1\le k \le m-1$ and $(\te_m\;
z_{m+r}) \red \sou{\m}. (z_j \,\bar{y})]$, or $[\exists j: t \red
\,\sou{\m}.(z_j\;\bar{y})]\}$, take also ${\cal{R}}= \{\bar{y}\}
\leadsto \cal{S}$.\\

%It is important to clarify the case $m =0$ in the definition of $\cal{S}$, this corresponds exaclty to $\exists\, j$: $t \, \red\, \sou{\m}.(z_j\;\bar{y})$, thus they do not exist terms $\te_i$.\\

It is clear that $\cal{S}$ is a $\mu$-saturated set. Let
${\mathcal{M}} = \langle {\cal{A}}, {\cal{S}},{\cal{R}} \rangle$ and
an $\mathcal{M}$-interpretation $I$ such that $I(X)= {\cal{R}}$. By
the corollary \ref{closed}, $E \in [({\cal{R}} \leadsto {\cal{S}})
\leadsto {\cal{R}}] \leadsto (\{\bar{y}\} \leadsto \cal{S})$. Let us
check that $x \in ({\cal{R}} \leadsto {\cal{S}}) \leadsto
{\cal{R}}$. For this, we take $\te \in ({\cal{R}} \leadsto \cal{S})$
and we prove that $(x\; \te) \in {\cal{R}}$, i.e, $((x\;
\te)\;\bar{y}) \in \cal{S}$.  By the definition of $\cal{S}$,
$(z_r\;\bar{y})\in \cal{S}$ for each $r \ge 0$, hence $z_r \in
{\cal{R}}$. Therefore $(\te\; z_r) \in \cal{S}$, so we have $\forall
r' \ge 0$:

\begin{enumerate}
\item Either $\,\exists m\ge 1, \,\exists \te_1,...,\te_m,\, \exists j:$ 
\begin{itemize}
\item $(\te\; z_r)\,\red\, {\sou{\m}}.((x\; \te_1)\;\bar{y})$
\item $(\te_k\; z_{k+r'})\,\red\, {\sou{\m}}.((x\; \te_{k+1})\;\bar{y})$ $\quad$ for every $1\le k \le m-1$
\item $(\te_m\; z_{m+r'})\,\red\,{\sou{\m}}.(z_j\,\bar{y})$.
\end{itemize}

More generally, since this holds for any $r'$, take $r'=r+1$, then,

$\exists m\ge 1, \,\exists \te_1,...,\te_m,\, \exists j:$ 
\begin{itemize}
\item $(\te\; z_r)\,\red\, {\sou{\m}}.((x\; \te_1)\;\bar{y})$
\item $(\te_k\; z_{k+1+r})\,\red\, {\sou{\m}}.((x\; \te_{k+1})\;\bar{y})$ $\quad$for every $1\le k \le m-1$
\item $(\te_m\; z_{m+1+r})\,\red\,{\sou{\m}}.(z_j\,\bar{y})$.
\end{itemize}

Therefore take $m'=m+1$, and the terms $\te_1'=\te ,\,\te_2'=\te_1,\,...,\te_{m+1}'=\te_m$, hence check easily that we have for any fixed $r$:

$\exists m'\ge 1, \,\exists \te_1',...,\te_{m'}',\, \exists j:$ 
\begin{itemize}
\item $((x\; \te)\;\bar{y})\,\red\,{\sou{\m}}.((x\; \te_1')\;\bar{y})$
\item $(\te_1'\; z_r)\,\red\, {\sou{\m}}.((x\; \te_2')\;\bar{y})$
\item $(\te_k'\; z_{k+r})\,\red\, {\sou{\m}}.((x\; \te_{k+1}')\;\bar{y})$ $\quad$ for every $1\le k \le m'-1$
\item $(\te_{m'}'\; z_{m'+r})\,\red\,{\sou{\m}}.(z_j\,\bar{y})$.
\end{itemize}

\item Or $\exists j: (\te\; z_r)\,\red\,{\sou{\m}}.(z_j\, \bar{y})$,
 then $((x\;\te)\;\bar{y})\, \red\, {\sou{\m}}.((x\;
 \te_1')\;\bar{y})$ and $(\te_1'\; z_r)\,\red\, {\sou{\m}}.(z_j\,
 \bar{y})$ with $m'=1$ and $\te'_1 =\te$. Therefore $((x\;
 \te)\;\bar{y}) \in \cal{S}$$)$.
\end{enumerate}

Thus $((x\; \te)\;\bar{y}) \in \cal{S}$ which implies that
$((E\,x)\,\bar{y}) \in \cal{S}$. By the fact that $E$ is a closed
term, the $\l$-variable $x$ and the sequence $\bar{y}$ are different
from each $z_i$, one can ensure that the assertion $[\exists j:
((E\,x)\,\bar{y})\, \red\, {\sou{\m}}.(z_j\, \bar{y})]$ can not
hold. Then for $r=0$, $\exists m\ge 1, \exists \te_1,...,\te_m,\exists
j$ such that:
\begin{itemize}
\item $ ((E\,x)\,\bar{y}) \,\red\,\sou{\m}. ((x \;\te_1)\,\bar{y})$
\item $(\te_k\; z_{k})\, \red\, \sou{\m}.((x \,\te_{k+1})\,\bar{y})$ $\quad$for every $1\le k \le m-1$ 
\item $(\te_m\; z_{m}) \,\red \, \sou{\m}. (z_j \,\bar{y})$ $\quad$ for some $1\le j\le m$.
\end{itemize}

\underline{Syntactical proof:}

Now we give a syntactical proof of this result. Let
$\bar{O}=O_1,...,O_n$ be new constants, $A= O_1,...,O_n \to \bot$ and
$\bar{y} = y_1...y_n$ a sequence of variables. By the lemma \ref{toto}
$\v_{\bar{O}} E:(\neg A \to A) \to A$, then, $x:\neg A \to A,
(y_i:O_i)_{1 \le i \le n} \v_{\bar{O}} (E\;x)\bar{y}:
\bot$. Therefore, $(T\;x)\bar{y} \red \t$, where $\t$ is a normal term
and $x:\neg A \to A, (y_i:O_i)_{1 \le i \le n} \v_{\bar{O}} \t:\bot$.

Following the form of $\t$ we have only one case to examine, the
others give always contradictions. This case is $\t=\sou{\m}.(x
\;U_1)\;t_1...t_n$ where $U_1,t_1,...,t_n$ are normal terms,
$x:\neg A \to A, (y_i:O_i)_{1 \le i\le n} \v_{\bar{O}} U_1:\neg
A\;;(b_j: \bot)_{1\le i\le m}$ and for all $1 \leq k \leq n$, $x:\neg
A \to A, (y_i:O_i)_{1 \le i \le n} \v_{\bar{O}}t_k:O_k\;;(b_j:
\bot)_{1\le j\le m}$. 
We deduce, by the lemma \ref{obstacle}, that, for all $1 \leq k \leq n$, $t_k = y_k$.

We prove, by induction and using the lemma \ref{obstacle}, that if
$x:\neg A \to A, (y_i:O_i)_{1 \le i \le n}, (z_k:A)_{1\leq k \leq i-1}
\v_{\bar{O}} U_i:\neg A\;;(b_j: \bot)_{1\le j\le m}$, then
\begin{displaymath}
\left\{\begin{array}{l} (U_i\; z_i) \red \sou{\m}.(x\; U_{i+1})\bar{y} \;\;and \; x:\neg A \to A, (y_i:O_i)_{1
\le i \le n}, (z_k:A)_{1\leq k \leq i} \v_{\bar{O}} U_{i+1}:\\ \neg A\;;(b_j: \bot)_{1\le j\le m}\\ \;or\;\\ \exists j : (1\le j \le i),\; such\;\; that:\; (U_i \;z_i) \red \sou{\m}.z_j\bar{y}
\end{array}\right.
\end{displaymath} 
The sequence $(U_i)_{i \geq 1}$ is not infinite, else the term $((E
\;\l x .\m a.(x\,z)) \bar{y})$ is not normalizable, which is
impossible, since\\ $x:\neg A,z:A, (y_i:O_i)_{1 \le i \le n}
\v_{\bar{O}} ((E \;\l x.\m a.(x\,z)) \bar{y}): \bot$.
\end{pro}

\begin{cor}\label{call/cc}
Let $E$ be a closed term of type $(\neg X \to X) \to X$, then, for
each term $u$, for each sequence $\bar{w} \in\mathcal{T}^{<\o}$ and
for each sequence $(v_i)_{i \in {\mathbb{N}}^*}$ of terms. There exist
$m \in \mathbb{N}$ and terms $\te_1,...,\te_m$ such that we have:
\begin{itemize}
\item$(E\; u) \bar{w} \red \sou{\mu}. (u\,\, \te_1)\;\bar{w} $
\item$(\te_i\; v_i)\red \sou{\mu}. (u\,\,
  \te_{i+1})\;\bar{w}$ $\;\;\;\;$ for all $ 1\le i \le m-1$
\item$(\te_m\; v_m) \red \sou{\mu}. (v_i\,\, \bar{w})$
$\;\;\;\;\;\;\;$ for some $ 1\le i \le m$
\end{itemize}
\end{cor}

\begin{pro} 
Immediate from the previous theorem and the lemma \ref{sigma}.
\end{pro}

\begin{re} 
In the {\tt C} programming language, there exist  ``escape''
instructions which allow to manage errors without stopping the
program. These are {\tt setjmp} and {\tt longjmp}. If we reduce $(E_1
\, \lambda y.h)\, \overline{w}$, we obtain $\mu
a.(a\,(h[y:=\theta_1]\,\overline{w}))$. When $\theta$ is executed with
some value $v$, the environment is restored and we get $(a.(v \,
\overline{w}))$. In other words, in the term $(E_1 \, \lambda y.h)$,
$E_1$ plays the role of the {\tt setjmp} instruction and occurences
of the variables $y$ in $h$ are the {\tt longjmp} instruction. The
corollary \ref{call/cc} says that every term of type $(\neg X \f X) \f
X$ has the same operational behavior of $E_1$ but often in
several steps $($the sequence of $\theta_i$$)$.
\end{re}

\section{Future work} 

Through this work, we have seen that the propositional types of the
system $S_\m$ are complete for the semantics defined previously.

\begin{enumerate}

\item What about the types of the second order typed $\l\m$-calculus?
We know that, for the system ${\cal F}$, the $\forall^+$-types (types
with positive quantifiers) are complete for a realizability semantics
(see \cite{4fn} and \cite{4Labib}). But for the classical system
${\cal{F}}$, we cannot generalize this result. We check easily that,
if $t=\m a.(a \, \l y_1 \l z \m b.(a \, \l y_2 \l x.z))$ and $A =
\forall \, Y \{Y \to \forall X (X \to X)\}$, then $t \in |A|$, but $t$
does not have the type $A$. This is due to the presence of $\forall$
in right-hand-side of $\to$, hence, we need to add more restrictions
on the positions of $\forall$ in the $\forall^+$-types to obtain a
smallest class of type that we suppose can be proved complete.

\item The problem is not the same when we consider the propositional
classical natural deduction system with the connectives $\wedge$ and
$\vee$. In previous works \cite{4NS2} and \cite{4NS1}, we define
interpretations of $\wedge$ and $\vee$ according to the functional
constructors $\curlywedge$ and $\curlyvee$ respectively as follows:

\begin{itemize}
\item ${\cal{K}} \curlywedge {\cal{L}} =\{t \in {\cal{T}}\,/\,
(t\,\pi_1) \in {\cal{K}}$ and $(t\,\pi_2) \in {\cal{L}} \}$

\item ${\cal{K}} \curlyvee {\cal{L}} =\{t \in {\cal{T}}\,/\,$ for each
$u, v$ if $($for each $r \in {\cal{K}}$, $s \in {\cal{L}}$ : $u[x:=r]
\in {\cal{S}}$ and $v[y:=s] \in {\cal{S}}$$)$, then $(t\,[x.u,y.v])
\in {\cal{S}}$$\}$
%where ${\cal{K}}$, ${\cal{L}}$ are sets of terms and ${\cal{S}}$ is a particular set of terms.
\end{itemize}

These interpretations allow to obtain a correctness result. We can
easily check that the term $\m a.(a\; \<\m b.(a\;\<\l x.x,\m c. (b\;\l
y.\l z.z)\>), \l x.x\> )$ belongs to the interpretation of the type $A
= (X \to X) \wedge (X \to X)$ but it does not have the type $A$. The
treatment of the disjunction is even a delicate matter, so we think that to
circumventing this difficulties, and if we hope a completeness
theorem, some deep modifications should be brought to our semantics.

\end{enumerate}
\bigskip

{\bf \large{Acknowledgements}}: We wish to thank R. Matthes and P. De
Groote for helpful discussions.

%%%%%%%%%%%%%%%%%%%%%%%%%%%%%%%%%%%%%%%%%%%%%%%%%%%
\section{Appendix}

This part is devoted to the proof of the lemma \ref{perdre patience}.

\begin{nt} 
Let $y$ be a $\l$-variable. The expression $u \rd_{\b y} v$ $($resp $u
\rd_{\m y}v $$)$ means
that we reduce in $u$ only a $\b$ $($resp $\m$$)$-redex where $y$ is the
argument, i.e, a redex in the form $(\l z.u\; y)$ $($resp $(\m b. u\;
y)$$)$. We denote by $\rd_y$ the union of $ \rd_{\b y}$ and $\rd_{\m y}$
and $ \red_y$ $($resp $\red_{\b y}$, $\red_{\m y}$$)$ the transitive and
reflexive closure of $ \rd_y$ $($resp $ \rd_{\b y}$, $\rd_{\m y}$$)$.
\end{nt}

\begin{lm}\label{y} 
Let $t$ be a normal term, $\s = [(a_i:=^*y)_{1\le i \le n}]$ and $\t$
the normal form of $t\s$, then, $t\s \red_y \t$.
\end{lm}

\begin{pro} 
By induction on the normal term $t$, the important case is the one
where $t= (a_i\;u)$ and $u$ a normal term, the others are direct
consequences of induction hypothesis.  Let us examine the different
forms of the normal term $u$, here there are two important subcases
$u= \l x.v$ and $u =\m b.v$ with $v$ a normal term $($these are the
two cases where there is creation of redexes after 
substitution$)$.
\begin{enumerate}
\item If $u= \l x.v$, then, $u\s= \l x.v\s$ and $t\s = (a_i\;(\l
x.v\s\;\;y))\, \rd_{\b y} (a_i\;v\{\s+[x:=y]\})$. By induction
hypothesis, $v\s \red_y v'$ where $v'$ is the normal form of $v\s$,
hence $(a_i\;v\{\s+[x:=y]\}) \,\red_y (a_i\;v'[x:=y])$ which is the normal
form of $t\s$.
\item If $u =\m b.v$, then, $u\s= \m b.v\s$ and $t\s= (a_i\;(\m
b.v\s\;\;y)) \,\rd_{\m y} (a_i\;\m b.v\{\s+[b:=^*y]\})$. By induction
hypothesis, $v\{\s+[b:=^*y]\}$ is normalizable only with $\red_y$
reductions, therefore $t\s$ is also normalizable only by $\red_y$
reductions.
\end{enumerate}
\end{pro}

\begin{lm}\label{yredex} 
Let $t$ be a normal term, $\t$ the normal form of $t[a:=^*y]$ and $A,B$ two types. If $\Gamma\, , y:A \vdash \t:B;
\Delta$. Then $\Gamma\, , y:A \vdash t[a:=^*y]:B; \Delta$.
\end{lm}

\begin{pro} 
By induction on the length of the reduction $t[a:=^*y]\red_y
\t$. By the lemma \ref{y}, it suffices to prove the  following
lemma.
\end{pro}

\begin{lm} 
Let $\t$ be a normal term, $t$ a term and $A,B$ two types. If $t
\rd_{\b y} \t$ $($resp $t \rd_{\m y} \t$$)$ and $\G\,, y:A \v
\t:B\;;\D$ then $\G\,, y:A \v t:B\;;\D$.
\end{lm}

\begin{pro} 
By induction on $t$, we examine how $t \rd_{\b y} \t$ (resp $t \rd_{\m
y} \t$). The proof is similar to the proof of (2) of the lemma
\ref{fatiguant}.
\end{pro}

\begin{lm}\label{coller_y} 
Let $t$ be a normal term, $y$ a $\l$-variable such that $y \not \in
Fv(t)$, $\s =[a:=^*y]$ and $A,B,C$ types. If $\Gamma\, , y:A\v t\s
:B\; ; \Delta, a: C$, then, $\Gamma\, \v t:B \; ; \Delta, a:A \to C$.
\end{lm}

\begin{pro} 
By induction on $t$.
\begin{enumerate}

\item $t =(x\,u_1)\,u_2...u_n$, then, $t\s =(x\,u_1\s)\,u_2\s...u_n\s$
and $\Gamma\,, y:A\v (x\,u_1\s)\,u_2\s...$\\$u_n\s:B\; ; \Delta, a:
C$. Therefore $x : E_1,...,E_n \to B \in \G$ and $\G, y:A \v u_i\s:E_i;
\Delta, a:C$. By induction hypothesis, we have  $\G \v u_i\s:E_i;
\Delta, a:A \to C$, hence $\Gamma\, \v (x\;u_1)\,u_2...u_n:B\; ; \Delta,
a: A \to C$.

\item $t= \l x.u$, then, $t\s= \l x.u\s$ and $\G\,, y:A \v \l x.u\s: B;
\Delta, a:C$, this implies that $B=F \to G$ and $\G\,, y:A , x: F\v u\s: G;
\Delta, a:C$. By induction hypothesis, $\G\,,x:F \v u:G; \Delta,
a: A \to C$, then, $\G\,\v \l x.u:F \to G; \Delta, a: A \to C$,
therefore $\G\,\v \l x.u:B; \Delta, a: A \to C$.

\item $t = \m b.u $, then, $t\s= \m b.u\s$ and
$\Gamma\, , y:A \v \m b.u\s :B\; ; \Delta, a: C$, this implies that
$\Gamma , y:A\v u\s:\perp\; ; \Delta, a: C, b:B$. By induction
hypothesis, $\Gamma\, \v u :\perp \; ; \Delta, a: A \to C, b:B$,
therefore $\Gamma\, \v \m b.u:B\; ; \Delta, a: A \to C$.

\item $t=(a\;u)$, then $t\s =(a\;(u\s\;y))$ and $\Gamma\, , y:A
\v (a\;(u\s\;y)) :\perp \; ; \Delta, a: C$, this implies that
$\Gamma\, , y:A \v (u\s\;y) : C \; ; \Delta, a: C$ and $\Gamma\, , y:A
\v u\s :A \to C \; ; \Delta, a: C$. By induction hypothesis,
$\Gamma\, \v u :A \to C \; ; \Delta, a: A \to C$, therefore $\Gamma\,
\v (a\;u) :\perp \; ; \Delta, a: A \to C$.

\item $t=(b\;u)$, then, $t\s= (b\;u\s)$ and $\Gamma\, , y:A \v (b\;u\s)
:\perp \; ; \Delta, a: C$, this implies that $\Gamma\, , y:A \v u\s :G
\; ; \Delta,b:G, a: C$. By induction hypothesis, $\Gamma\, \v u :G
\; ; \Delta,b:G, a: A \to C$, therefore $\Gamma\, \v (b\;u) :\perp \;
; \Delta, a: C$.

\end{enumerate}
\end{pro}

\begin{pro}[of lemma \ref{perdre patience}]  
By induction on $t$, the cases where $t = (x\;u_1)\,u_2...u_n$ and $t
= \l x.u$ are similar to those in the proof of (2) of the lemma
\ref{fatiguant}. Let us examine the case where $t=\m a.u$, then $(t\;
y) \red \m a.u[a:=^*y] \red \m a.u' = \t$ where $u'$ is the normal
form of $u[a:=^*y]$. We have $\Gamma\, , y:A \v \m a.u':B; \;
\Delta$, then $\Gamma\, , y:A \v u': \perp; \;
\Delta, a:B$. By the lemma \ref{y}, $ u[a:=^*y] \red_{y} u'$,
then, by the lemma \ref{yredex}, $\Gamma\, , y:A \v u[a:=^*y]:
\perp; \; \Delta, a:B$. Hence by the lemma \ref{coller_y},
$\Gamma\, \v u: \perp; \; \Delta, a:A \to B$ finally $\Gamma\, \v
\m a.u:A\to B; \; \Delta$.
\end{pro}


\begin{thebibliography}{99}



\bibitem{4coq} T. Coquand {\em Completeness theorem and $\l$-calculus.}
The 7th International Conference, TLCA 2005, Nara, Japan, April 21-23,
2005, pp. 1-9, volume 3461/2005.



\bibitem{4davnou2} R. David and K. Nour. {\em A short proof of the
strong normalization of the simply typed $\l \m$-calculus.}
Schedae Informaticae vol 12, pp. 27-33, 2003.

\bibitem{4dav} R. David. {\em Une preuve simple de r\'esultats classiques en $\l$-calcul.} Compte Rendu de l'Acad\'emie des Sciences. Paris, Tome 320, S\'erie 1, pp. 1401-1406, 1995.


\bibitem{4fn} S. Farkh and K. Nour. {\em Un r\'esultat de compl\'etude
  pour les types $\forall^+$ du syst\`eme $\cal{F}$.} CRAS. Paris 326, S\'erie I, pp. 275-279, 1998.

\bibitem{4fn2} S. Farkh and K. Nour. {\em Types Complets dans une
  extension du syst\`eme ${\cal{AF}}2$.} Informatique Th\'eorique et
  Application, 31-6, pp. 513-537, 1998.

\bibitem{4Gir} J.-Y. Girard, Y. Lafont, P. Taylor.  Proofs and types.
{\it Cambridge University Press}, 1986.

\bibitem{4Grif} T. Griffin. {\em A formulae-as-types notion of
control.} Proc. POLP, 1990.  

\bibitem{4hind1} J. R. Hindley. {\em The simple semantics for
  Coppe-Dezani-Sall\'e types.} Proceeding of the 5th Colloquium on
  International Symposium on Programming, pp. 212-226, April 06-08, 1982.

\bibitem{4hind2} J. R. Hindley. {\em The completeness theorem for typing
  $\l$-terms.} Theoretical Computer Science, 22(1), pp. 1-17, 1983.

\bibitem{4hind3} J. R. Hindley. {\em Curry's type-rules are complete
  with respect to the F-semantics too.} Theoretical Computer Science,
  22, pp. 127-133, 1983.

\bibitem{4kam} F. Kamareddine and K. Nour. {\em A completeness result for
  a realizability semantics for an intersection type system.} Annals of Pure and Applied Logic, vol 146, pp. 180-198, 2007

\bibitem{4kri1} J.-L. Krivine. {\em Lambda calcul, types et mod\`eles.} Masson, Paris, 1990.

\bibitem{4kri4} J.-L. Krivine. {\em Op\'erateurs de mise en m\'emoire
  et traduction de G\"odel.} Archive for Mathematical Logic, vol 30,
  pp. 241-267, 1990.

\bibitem{4Labib} R. Labib-Sami. {\em Typer avec (ou sans) types
  auxiliaires.} Manuscrit, 1986.

\bibitem{4NS2} K. Nour and K. Saber. {\em A Semantics of Realizability for the Classical Propositional Natural Deduction.} Electronic Notes in
  Theoretical Computer Science, vol 140, pp. 31-39, 2005.

\bibitem{4NS1} K. Nour and K. Saber. {\em A semantical proof of strong
  normalization theorem for full propositional classical natural
  deduction.} Archive for Mathematical Logic, vol 45, pp. 357-364, 2005.

\bibitem{4Nour5} K. Nour. {\em Op\'erateurs de mise en m\'emoire et
  types $\forall$-positifs.} Theoretical Informatics and Applications,
  vol 30, n° 3, pp. 261-293, 1996.

\bibitem{4Nour6} K. Nour. {\em Mixed Logic and Storage Operators.} Archive for Mathematical Logic, vol 39, pp. 261-280, 2000.

\bibitem{4Par1} M. Parigot {\em $\l \m$-calculus: An algorithm
  interpretation of classical natural deduction.}  Lecture Notes in
  Artificial Intelligence, vol 624,
pp. 190-201. Springer Verlag, 1992.

\bibitem{4Par2} M. Parigot. {\em Proofs of strong normalization for
second order classical natural deduction.} Journal of Symbolic
Logic, vol 62 (4), pp.  1461-1479, 1997.

\bibitem{4Py} W. Py. {\em Confluence en $\l \m$-calcul.} PhD thesis, University of Chamb\'ery, 1998.

\bibitem{Moi} K. Saber. {\em \'Etude d'un $\l$-calcul issu d'une logique classique} PhD Thesis, University of Chamb\'ery, 2007.

\bibitem{4Sau} A. Saurin. {\em Separation and the
  $\l\m$-calculus.} Proceedings of the Twentieth Annual IEEE Symp. on
Logic in Computer Science, {LICS} 2005, IEEE Computer Society Press, pp. 356-365, 2005.


\bibitem{4Tait} W. W. Tait, {\em A realizability interpretation of the
  theory of species.} In : R. Parikh (Ed.), Logic Colloquium Boston
  1971/72, vol. 435 of Lecture Notes in Mathematics, Springer Verlag, pp. 240-251, 1975.


\end{thebibliography}
\end{document}